\documentclass{article}
\usepackage{amssymb,latexsym}
\usepackage{amsmath}
\usepackage{graphics}
\usepackage{graphicx}
\newtheorem{theorem}{Theorem}

\numberwithin{equation}{section} \numberwithin{theorem}{section}
\numberwithin{corollary}{section} \numberwithin{lemma}{section}
\numberwithin{remark}{section} \numberwithin{definition}{section}

\begin{document}
\title{On Timelike Tubular Weingarten Surfaces in Minkowski 3-Space}
\author{A. Z.  Azak$^{a}$, M. Masal$^{a}$, S. Hal{\i}c{\i}$^{b}$}

\date{}

\maketitle
\begin{center}
$^{a}$ Department of Mathematics Teaching , Faculty of Education \\
Sakarya University,  Sakarya/TURKEY

$^{b}$ Department of Mathematics, Faculty of Arts and Sciences \\
Sakarya University, Sakarya/TURKEY

\end{center}

\begin{abstract}
In this paper, we study the timelike tubular Weingarten
surfaces in 3-dimensional Minkowski space $IR_1^3 $.We have
obtained some conditions for being $\left( {K_{II} ,H} \right)$,
$\left( {K_{II} ,K} \right)$, timelike tubular Weingarten surfaces
where   are the second Gaussian curvature the Gaussian curvature
and the mean curvature, respectively.
\end{abstract}

\section{Introduction}\label{S:intro}

TLet $f$  and  $g$ be smooth functions on a surface $M$ in
3-dimensional Minkowski space $IR_1^3 $. The existence of a
nontrivial functional relation $\Phi \left( {f,g} \right) = 0$
namely the Jacobian determinant of $f$ and $g$ functions $\Phi
\left( {f,g} \right) = \det \left( {\begin{array}{*{20}c}
   {f_s } & {f_t }  \\
   {g_s } & {g_t }  \\
\end{array}} \right) = 0$,where $f_s  = \frac{{\partial f}}{{\partial s}},\;f_t  = \frac{{\partial f}}{{\partial
t}}$. A surface is called a Weingarten surface, if there is a
nontrivial relation $\Phi \left( {k_1 ,k_2 } \right) = 0$ between
its principal curvatures $k_1 $ and $k_2 $ , or, equivalently
there is a nontrivial relation $\Phi \left( {K,H} \right) = 0$
between its Gaussian curvature $K$ and mean curvature $H$. Also we
call a surface as a linear Weingarten surface such that the linear
combination $aK + bH = c$ is constant along each ruling, where
$a,b,c \in ,\;\left( {a,b,c} \right) \ne \left( {0,0,0} \right)$.
Several geometers have studied Weingarten and linear Weingarten
surfaces. For non-developable surfaces Kühnel studied $\left\{
{H,K_{II} } \right\}$ and $\left\{ {K,K_{II} } \right\}$
Weingarten surfaces \cite{Ku}. Then, G. Stamou extended this
article and gave linear Weingarten surfaces which satisfy $aK_{II}
+ bH + cH_{II}  = d$ is constant along rulings where $a^2  + b^2 +
c^2  \ne 0$ \cite{St}. Also in 3-dimensional Minkowski space
linear Weingarten surfaces which is foliated by pieces of circles
and linear Weingarten helicoidal surfaces under cubic screw motion
studied in \cite{Fe}, \cite{Ka}. F. Dillen and W. Sodsiri examined
$\left\{ {K,K_{II} ,H,H_{II} } \right\}$ Weingarten and linear
Weingarten surfaces in 2005 for the 3-dimensional Minkowski space
\cite{Di1}, \cite{Di2}, \cite{Di3}.

A tubular Weingarten and linear Weingarten surface were studied by
Ro and Yoon in 3-dimensional Euclidean space $E^3 $ \cite{Ro}.
Karacan and Bukcu constructed tubular surfaces with the assistance
of an alternative moving frame \cite{Kr2}. Then geodesics and
singular points of tubular surfaces are researched by using one
parameter spatial motion along a curve in Minkowski 3-space
\cite{Kr1}, \cite{Kr3}.

In this paper, we give some theorems and conclusions related to
timelike tubular Weingarten and linear Weingarten surfaces in
3-dimensional Minkowski space $IR_1^3 $.

\section{Preliminaries}\label{S:intro}

The Minkowski 3-space $IR_1^3 $ is the Euclidean 3-space $IR^3 $
provided with the indefinite inner product given by
$$\left\langle , \right\rangle  =  - dy_1^2  + dy_2^2  +
dy_3^2,$$ where $\left( {y_1 ,y_2 ,y_3 } \right)$ is natural
coordinates of $IR_1^3 $. Since $\left\langle , \right\rangle $ is
indefinite inner product, recall that a vector $\beta  \in IR_1^3$
can have one of the three causal characters it can be spacelike if
$\left\langle {\beta ,\beta } \right\rangle  > 0$ or $\beta  = 0$,
timelike if $\left\langle {\beta ,\beta } \right\rangle  < 0$ and
null (lightlike) if $\left\langle {\beta ,\beta } \right\rangle  =
0$ and $\beta  \ne 0$. Similarly, an arbitrary curve $\gamma  =
\gamma \left( t \right)$ in $IR_1^3 $ can locally called as
timelike, if its velocity vector $\gamma '\left( t \right)$ is
timelike. Recall that the norm of a vector is given by $\left\|
\beta  \right\| = \sqrt {\left| {\left\langle {\beta ,\beta }
\right\rangle } \right|} $ and that the timelike $\gamma \left( t
\right)$ is said to be of unit speed if $\left\langle {\gamma
'\left( t \right),\gamma '\left( t \right)} \right\rangle  =  -
1$. Morever, the velocity of curve $\gamma \left( t \right)$ is
the function $\upsilon \left( s \right) = \left\| {\gamma '\left(
t \right)} \right\|$. Denote by $\left\{ {t,n,b} \right\}$ the
moving Frenet frame along the curve $\gamma \left( t \right)$ in
the Minkowski space $IR_1^3 $. Then Frenet formula of $\gamma
\left( t \right)$ in the space $IR_1^3 $ is defined by \cite{Pe}.
$$\begin{array}{l}
 t' = \kappa n \\
 n' = \kappa t + \tau b \\
 b' =  - \tau n \\
 \end{array}$$
where the prime denotes the differentiation with respect to $t$
and we denote by $\kappa ,\tau $ the curvature and the torsion of
the curve $\gamma $. Since $\gamma $ is a timelike curve
$$\begin{array}{l}
 \left\langle {t,t} \right\rangle  =  - 1,\left\langle {b,b} \right\rangle  = \left\langle {n,n} \right\rangle  = 1, \\
 \left\langle {t,n} \right\rangle  = \left\langle {t,b} \right\rangle  = \left\langle {b,n} \right\rangle  = 0. \\
 \end{array}$$
The vector product of the vectors $\beta  = \left( {\beta _1
,\beta _2 ,\beta _3 } \right)$ and $\mu  = \left( {\mu _1 ,\mu _2
,\mu _3 } \right)$ is defined by
$$\beta  \wedge \mu  = \left( {\beta _3 \mu _2  - \beta _2 \mu _3 ,\beta _3 \mu _1  - \beta _1 \mu _3 ,\beta _1 \mu _2  - \beta _2 \mu _1 } \right).$$
We denote a timelike surface $M$ in $IR_1^3 $ by
$$x\left( {s,t} \right) = \left( {x_1 \left( {s,t} \right),x_2 \left( {s,t} \right),x_3 \left( {s,t} \right)} \right)$$

Let $U$ be the standard unit normal spacelike vector field on a
surface $M$ defined by $U = \frac{{x_s  \wedge x_t }}{{\left\|
{x_s  \wedge x_t } \right\|}}$, where $x_s  = \frac{{\partial
x\left( {s,t} \right)}}{{\partial s}}$ and $x_t  = \frac{{\partial
x\left( {s,t} \right)}}{{\partial t}}$. Then the first fundamental
form $I$ and the second fundamental form $II$ of a timelike
surface  $M$ are defined by, respectively
$$I = Eds^2  + 2Fdsdt + Gdt^2, $$
$$II = eds^2  + 2fdsdt + gdt^2, $$
where
$$E = \left\langle {x_s ,x_s } \right\rangle ,\quad \;F = \left\langle {x_s ,x_t } \right\rangle ,\quad \;G = \left\langle {x_t ,x_t } \right\rangle $$
$$e = \left\langle {x_{ss} ,U} \right\rangle ,\;\quad f = \left\langle {x_{st} ,U} \right\rangle ,\;\quad g = \left\langle {x_{tt} ,U} \right\rangle. $$
On the other hand, the Gaussian curvature $K$ and the mean
curvature $H$ are given by, respectively
$$\begin{array}{l}
 K = \frac{{eg - f^2 }}{{EG - F^2 }}\left\langle {U,U} \right\rangle  \\
 H = \frac{{eG - 2fF + gE}}{{2\left( {EG - F^2 } \right)}}\left\langle {U,U} \right\rangle . \\
 \end{array}$$

From Brioschi's formula in a Minkowski 3-space, we can compute
$K_{II} $ of a surface by replacing the components of the first
fundamental form $E,\,F,\,G$ by the components of the second
fundamental form $e,\,f,\,g$ respectively in Brioschi's formula.
Consequently, the second Gaussian curvature $K_{II} $ of a
non-developable surface is defined by \cite{So}
$$ K_{II}  =
\frac{1}{{\left( {eg - f^2 } \right)^2 }}\left\{ {\left|
{\begin{array}{*{20}c}
   { - \frac{1}{2}e_{tt}  + f_{st}  - \frac{1}{2}g_{ss} } & {\frac{1}{2}e_s } & {f_s  - \frac{1}{2}e_t }  \\
   {f_t  - \frac{1}{2}g_s } & e & f  \\
   {\frac{1}{2}g_t } & f & g  \\
\end{array}} \right| - \left| {\begin{array}{*{20}c}
   0 & {\frac{1}{2}e_t } & {\frac{1}{2}g_s }  \\
   {\frac{1}{2}e_t } & e & f  \\
   {\frac{1}{2}g_s } & f & g  \\
\end{array}} \right|} \right\}.$$

\section{Timelike Tubular Surfaces of Weingarten Types
}\label{S:intro}

Let $\gamma :\left( {a,b} \right) \to IR_1^3 $ be a smooth unit
speed timelike curve of finite lenght which is topologically
imbedded in $IR_1^3 $. The total space $N_\gamma  $ of the normal
bundle of $\gamma \left( {\left( {a,b} \right)} \right)$ in
$IR_1^3 $ is naturally diffeomorphic to the direct product $\left(
{a,b} \right) \times IR_1^3 $ via the translation along $\gamma $
with respect to the induced normal connection. For sufficiently
small $r > 0$, the tubular surface of radius $r$ about the curve
$\gamma $ is the set:
$$T_r \left( \gamma  \right) = \left\{ {\left. {\exp _{\gamma \left( t \right)} \vartheta } \right|\;\vartheta  \in N_{_{\gamma \left( t \right)} } ,\;\left\| \vartheta  \right\| = r,\;a < t < b} \right\}.$$
For a sufficiently small the tubular timelike surface
$T_r\left(\gamma  \right) $ is a smooth surface in $IR_1^3 $. Then
the parametric equation of the timelike tubular surface
$T_r\left(\gamma  \right) $ can be expressed as
\begin{equation}\label{E2}
x\left( {t,\theta } \right) = \gamma \left( t \right) + r\left(
{\cos \theta {\kern 1pt} n + \sin \theta {\kern 1pt} b} \right)
\end{equation}
Furthermore, we have the natural frame
$\left\{{x_t,x_\theta}\right\}$ is given by
\begin{equation}\label{E3}
x_t  = \left( {1 + r\kappa \cos \theta } \right)t + r\tau \left(
{\cos \theta b - \sin \theta n} \right)
= \alpha t + r\tau v,\quad x_\theta   = r\left( {\cos \theta b - \sin \theta n} \right) = rv, \\
\end{equation}
where we put $\alpha  = 1 + r\kappa \cos \theta $ and
$v=\cos\theta \,b - \sin \theta \,n$. From which the components of
the first fundamental form are
\begin{equation}\label{E4}
E =  - \alpha ^2  + r^2 \tau ^2 ,\quad F = r^2 \tau ,\quad G =r^2.
\end{equation}
On the other hand, the unit normal spacelike vector field $U$ is
obtained by
$$U = \frac{{x_t  \wedge x_\theta  }}{{\left\| {x_t  \wedge x_\theta  } \right\|}} =  - \cos \theta \,n - \sin \theta \,b,$$
from this, the components of the second fundamental form of $x$
are given by
$$e = r\tau ^2  - \kappa \alpha \cos \theta ,\quad f = r\tau ,\quad g = r.$$
If the second fundamental form is non-degenerate,
$eg-f^2\ne0$,that is,$\kappa ,\alpha $ and $\cos \theta $ are
nowhere vanishing. In this case, we can define formally the second
Gaussian curvature $K_{II} $ on $T_r \left( \gamma  \right)$. On
the other hand, the Gauss curvature $K$,the mean curvature $H$ and
the second Gaussian curvature $K_{II} $ are given by, respectively
\begin{equation}\label{E5}
\begin{array}{l}
K = \frac{{\kappa \cos \theta }}{{r\alpha }},
\end{array}
\end{equation}
\begin{equation}\label{E6}
H = \frac{{ - \left( {1 + 2r\kappa \cos \theta }
\right)}}{{2r\alpha }},
\end{equation}
\begin{equation}\label{E7}
K_{II}  = \frac{1}{{4r\alpha ^2 \cos ^2 \theta }}\left( {4r^2
\kappa ^2 \cos ^4 \theta  + 6r\kappa \cos ^3 \theta  + \cos ^2
\theta  + 1} \right).
\end{equation}
Differentiating $K,H$ and $K_{II} $ with respect to $t$ and
$\theta $, we get
\begin{equation}\label{E8}
K_t  = \frac{{\kappa '\cos \theta }}{{r\alpha ^2 }},\quad \quad
K_\theta   =  - \frac{{\kappa \sin \theta }}{{r\alpha ^2 }},
\end{equation}
\begin{equation}\label{E9}
H_t  =  - \frac{{\kappa '\cos \theta }}{{2\alpha ^2 }},\quad \quad
H_\theta   = \frac{{\kappa \sin \theta }}{{2\alpha ^2 }},
\end{equation}

\begin{equation}\label{E9}
\begin{array}{l}
 \left( {K_{II} } \right)_t  = \frac{1}{{4r\alpha ^4 \cos \theta }}(2r^3 \kappa ^2 \kappa '\cos ^4 \theta  + 6r^2 \kappa \kappa '\cos ^3 \theta  + 4r\cos ^2 \theta  \\
 \,\,\,\,\,\,\,\,\,\,\,\,\,\,\,\,\,\,\,\,\,\,\,\,\,\,\,\,\,\,\,\,\,\,\,\,\,\,\,\,\,\,\,\,\,\,\,\,\,\,\,\, - 2r^2 \kappa \kappa '\cos \theta  - 2r\kappa '), \\
 \left( {K_{II} } \right)_\theta   = \frac{1}{{4r\alpha ^4 \cos ^4 \theta }}\left( { - 2r^3 \kappa ^3 \cos ^6 \theta \sin \theta  - 6r^2 \kappa ^2 \cos ^5 \theta \sin \theta } \right. \\
 \quad \quad \;\;\quad \quad \quad \quad \quad \quad\,\,\,\,\,\,\,\,\,  - 4r\kappa \cos ^4 \theta \sin \theta  + 4r^2 \kappa ^2 \cos ^3 \theta \sin \theta  \\
 \quad \quad \;\;\left. {\quad \quad \quad \quad \quad \quad\,\,\,\,\,\,\,\,\,  + 6r\kappa \cos ^2 \theta \sin \theta  + 2\cos \theta \sin \theta } \right). \\
 \end{array}
\end{equation}

Now, we investigate a tubular timelike surface
$T_r\left(\gamma\right) $ in $IR_1^3 $ satisfying the Jacobi
equation $\Phi \left( {X,Y} \right) = 0$. By using (3.7) and
(3.8), $\Phi \left( {X,Y} \right) = 0$ satisfies identically the
Jacobi equation $\Phi \left( {K,H} \right) = K_t H_\theta   -
K_\theta  H_t  = 0$. Therefore, $T_r\left(\gamma\right) $ is a
Weingarten surface. We consider a timelike tubular
$T_r\left(\gamma\right) $ with non-degenerate second fundamental
form in $IR_1^3 $ satisfying the Jacobi equation
\begin{equation}\label{E11}
\begin{array}{l}
\Phi \left( {K,K_{II} } \right) = K_t \left( {K_{II} }
\right)_\theta   - K_\theta  \left( {K_{II} } \right)_t  = 0
\end{array}
\end{equation}
with respect to the Gaussian curvature $K$ and the second Gaussian
curvature $K_{II} $. Then, by (3.7) and (3.9) equation (3.10)
becomes
$$r^2 \kappa ^2 \kappa '\cos ^2 \theta \sin \theta  + 2r\kappa \kappa '\cos \theta \sin \theta  + \kappa '\sin \theta  = 0$$
Since this polynomial is equal to zero for every $\theta $, all
its coefficients must be zero. Therefore, we conclude that
$\kappa'=0$. We suppose that a timelike tubular
$T_r\left(\gamma\right)$ with non-degenerate second fundamental
form in $IR_1^3 $ is $\left( {H,K_{II} } \right)$-Weingarten
surface. Then it satisfies the equation
\begin{equation}\label{E12}
\begin{array}{l}
H_t \left( {K_{II} } \right)_\theta   - H_\theta  \left( {K_{II} }
\right)_t  = 0,
\end{array}
\end{equation}
which implies
\begin{equation}\label{E13}
r^2 \kappa ^2 \kappa '\cos ^2 \theta \sin \theta  + 2r\kappa
\kappa '\cos \theta \sin \theta  + \kappa '\sin \theta  = 0
\end{equation}
from (3.12) we can obtain $\kappa ' = 0.$ \\Consequently, we have
the following theorems:
\begin{theorem} \label{T3.1.}A timelike tubular surface in a Minkowski 3-space is a Weingarten surface.
\end {theorem}
\begin{theorem} \label{T3.2.}Let $\left( {X,Y} \right) \in \left\{ {\left( {K,K_{II} } \right),\left( {H,K_{II} } \right)} \right\}$
and let $T_r \left( \gamma  \right)$ be a timelike tubular surface
in Minkowski 3-space with non-degenerate second fundamental form.
If $T_r \left( \gamma  \right)$ is a
$\left({X,Y}\right)$-Weingarten surface, then the curvature of
$T_r \left( \gamma  \right)$ is a non-zero constant.
\end {theorem}
Finally, we study a timelike tubular $T_r \left( \gamma  \right)$
in $IR_1^3 $ is a linear Weingarten surface, that is, it satisfies
the equation
\begin{equation}\label{E14}
\begin{array}{l}
aK + bH = c.
\end{array}
\end{equation}
Then, by (3.4) and (3.5) we have
$$\left( {2a\kappa  - 2br\kappa  - 2r^2 c\kappa } \right)\cos \theta  - b - 2rc = 0.$$
Since $\cos \theta $ and $1$ are linearly independent, we get
$$2a\kappa  - 2br\kappa  - 2r^2 c\kappa  = 0,\quad b =  - 2rc,$$
which imply
$$\kappa \left( {a + cr^2 } \right) = 0.$$
If $a + cr^2  \ne 0$,then $\kappa  = 0.$ Thus, $T_r \left( \gamma
\right)$ is an open part of a circular cylinder.\\
Next, suppose that a timelike tubular $T_r \left( \gamma \right)$
with non-degenerate second fundamental form in $IR_1^3 $ satisfies
the equation
\begin{equation}\label{E15}
aK + bK_{II}  = c.
\end{equation}
By (3.4) and (3.6), equation (3.14) becomes
$$\left( {4ar\kappa ^2
+ 4br^2 \kappa ^2  - 4cr^3 \kappa ^2 } \right)\cos ^4 \theta  +
\left( {4a\kappa  + 6br\kappa  - 8cr^2 \kappa } \right)\cos ^3
\theta  + \left( {b - 4cr} \right)\cos ^2 \theta  + b = 0.$$ Since
the identity holds for every $\theta $, all the coefficients must
be zero. Therefore, we have
$$\begin{array}{l}
 \quad \quad\ 4ar\kappa ^2  + 4br^2 \kappa ^2  - 4r^3 c\kappa ^2  = 0, \\
 \quad\quad\quad \quad\quad 4a\kappa  + 6br\kappa  - 8cr^2 \kappa  = 0, \\
 \quad\quad \quad\quad \quad\quad \quad \quad\quad \quad\ b - 4rc = 0, \\
\quad\quad \quad\quad \quad\quad \quad \quad \quad \quad \quad\quad\ b = 0. \\
 \end{array}$$
Thus, we get $b = 0,c = 0$ and $\kappa  = 0.$ In this case, the
second fundamental form of $T_r \left( \gamma  \right)$ is
degenerate. \\Suppose that a timelike tubular $T_r \left( \gamma
\right)$ with non-degenerate second fundamental form in $IR_1^3 $
satisfies the equation
$$aH + bK_{II}  = c.$$
By (3.5), (3.6) and (3.15), we have
$$\left( { - 4ar^2 \kappa ^2  + 4br^2 \kappa ^2  - 4cr^3 \kappa ^2 } \right)\cos ^4 \theta  + \left( { - 6ar\kappa  + 6br\kappa  - 8cr^2 \kappa } \right)\cos ^3 \theta  + \left( { - 2a + b - 4cr} \right)\cos ^2 \theta  + b = 0.$$
from which we can obtain $b = 0$ and $\kappa  = 0.$
\\Consequently, we have the following theorems:
\begin{theorem} \label{T3.3.}Let $T_r \left( \gamma  \right)$ be a timelike tubular surface satisfying the linear equation $aK + bH = c.$
If, $a + br \ne 0,$ then it is an open part of a circular
cylinder.
\end {theorem}
\begin{theorem} \label{T3.4.}Let $\left( {X,Y} \right) \in \left\{ {\left( {K,K_{II} } \right),\left( {H,K_{II} } \right)} \right\}.$
Then there are no $\left( {X,Y} \right)$-linear Weingarten tubular
in Minkowski 3-space $IR_1^3 $ .
\end {theorem}


\begin{thebibliography}{9}



\bibitem{Di1} F. Dillen, W. K\"{u}hnel, \emph{Ruled Weingarten Surfaces in Minkowski 3-space}, Manuscripta Math. 98, (1999), 307-320.


\bibitem{Di2} F. Dillen, W.Sodsiri,  \emph{Ruled Surfaces of Weingarten Type in Minkowski 3-space}, J. Geom. 83, (2005), 10-21.



\bibitem{Di3} F. Dillen, W.Sodsiri, \emph{Ruled surfaces of Weingarten type in Minkowski 3-space II}, J. Geom. 84, (2005), 37-44.


\bibitem{Fe} J. Fenghui , Y. Wang ,  \emph{Linear Weingarten Helicoidal Surfaces in Minkowski 3-space}, Differential Geometry-Dynamical
Systems 12,(2010),95-101.


\bibitem{Ka} \"{O}. B. Kalkan , R. Lopez , D. Saglam , \emph{Linear Weingarten Surfaces Foliated by Circles in Minkowski Space},
arxiv:0909.2552v1, (2009).


\bibitem{Kr1} M. K. Karacan , H. Es , Y. Yaylý  \emph{Singular Points of Tubular Surfaces in Minkowski 3-space}, Sarajevo Journal of
Mathematics (14) 2, (2006),73-82 .


\bibitem{Kr2} M. K. Karacan , B. Bukcu ,  \emph{An Alternative Moving Frame
for Tubular Surfaces Around Timelike Curves in Minkowski
3-space},Balkan Journal of Geometry and Its Applications, (2),
12,(2007), 73-80.

\bibitem{Kr3} M. K. Karacan , Y. Yaylý ,  \emph{On the Geodesics of Tubular Surfaces in Minkowski 3- space}, Bull. Malays. Math. Sci. Soc. (2)
31, (2008),1-10.


\bibitem{Ku} W. K\"{u}hnel ,  \emph{Ruled W-surfaces}, Arch. Math. 62,
(1994), 475-480 .


\bibitem{Pe} M. Petrovic-Torgasev, E. Sucurovic, \emph{Some Characterizations
of the Spacelike, the Timelike and the Null Curves on the
Pseudohyperbolic Space $H_0^2 $ in $E_1^3 $ }, Kragujevac J. Math.
22,(2000), 71-82.



\bibitem{Ro} J. S. Ro, D. W. Yoon,  \emph{Tubes of Weingarten Types in a
Euclidean 3-Space}, Journal of the Chungcheong Mathematical
Society 22, (2009),359-366.


\bibitem{So} W. Sodsiri ,   \emph{Ruled Surfaces of Weingarten Type in
Minkowski 3-space}, Doctorate thesis, Katholieke Universiteit
Leuven, 2005.
\bibitem{St} G. Stamou ,   \emph{Regelflachen vom Weingarten Typ}, Coloq.
Math. 79,(1999),77-84.

\end{thebibliography}
\end{document}